\documentclass[11pt]{article}

\usepackage{alltt}
\newcommand{\superscript}[1]{\ensuremath{^\textrm{#1}}}
\newcommand{\subscript}[1]{\ensuremath{_\textrm{#1}}}
\title{Recursion and Combinatorial Mathematics in \textit{Chandash\=astra}}

\author{Amba Kulkarni\\
Department of Sanskrit Studies,\\
University of Hyderabad\\
Hyderabad, India\\
apksh@uohyd.ernet.in\\}

\begin{document}
\maketitle
\begin{abstract}
Contribution of Indian Mathematics since Vedic Period has been recognised by the historians. \textit{Pi\.nga\d{l}a} (200 BC) in his book on '\textit{Chandash\=astra}', a text related to the description and analysis of meters in poetic work, describes algorithms which deal with the Combinatorial Mathematics.  These algorithms essentially deal with the binary number system -– counting using binary numbers, finding the value of a binary number,  finding the value of \superscript{n}C\subscript{r},  evaluating 2\superscript{n}, etc. All these algorithms are tail recursive in nature. Some of these algorithms also use the concept of  stack variables to stack the intermediate results for later use. Later work by \textit{\textit{Ked\=ar} Bhatt} (around 800 AD), however, has only iterative algorithms for the same problems. We describe both the recursive as well as iterative algorithms in this paper and also compare them with the modern works. 
\end{abstract}

\section{Introduction}
'Without any purpose, even a fool does not get initiated.' \footnote{\textit{prayojanam anuddishya na mando api pravartate}} Thus goes a saying in Sanskrit. If we look at the rich Sanskrit knowledge base, we find that all the branches of knowledge that exist in Sanskrit literature were originated in order to address some problems in day today life. While addressing them, we find that, there were also remarkable efforts in generalising the results and findings. For example, about the\textit{P\=a\d{n}ini}'s monumental work on \textit{asht\=adhy\=ayi} (500 BC), Paul Kipasky says "many of the insights of \textit{P\=a\d{n}ini}'s graamar still remain to be recaptured, but those that are already understood contribute a major theoretical contribution." (in the encyclopaedia of Language and Linguistics, ed Asher, pp 2923). 

Mathematics is also no exception to it. Contribution of Indian mathematicians dates back to the Vedic period\superscript{1}. The early traces of geometry and algebra are found in Sulbasutras of vedic period where the purpose of this geometrical and algebraic exercise was to build brick altars of different shapes to perform vedic rituals. Fixing Luni-Solar calendar was another important task which led to the development of calculus in India. The development florished in the classical period from Aryabhatta (500 AD) to Bhaskaracharya II (1150 AD) and further in the Kerala school of mathematics from 1350 AD to 1650 AD.

However the discovery of binary number system by Indians escaped the attention of Western scholars,  may be because  '\textit{Chandash\=astra}' was considered as mainly a text related to description and analysis of meters in poetic literary work, totally unrelated to mathematics. B. Van Nooten\superscript{2} brought it into limelight. It was the German philosopher Gotfried Leibniz in 1695 who re-discovered  the binary number system.

Binary number system forms the basis of the modern digital world and modern logic. The two numbers “0” and “1” may be interpreted as either “on” and “off”  or  “true” and “false", or “magnetised” and “demagnetised”, or “guru and laghu” as in the case of  '\textit{chandash\=astra}'  etc. 

      	\textit{Veda}s are in the form of poetic verses. Since the \textit{Veda}s have been transmitted orally from generations to generations, lot of care has been taken to preserve them in pure form, altogether avoiding any kind of destruction in their contents. Several methods have been adopted to preserve their contents. \textit{Veda}s are written in different “meters” (\textit{Chanda}s). These \textit{Chanda}s  have been studied in great detail. The \textit{Chandash\=astra} forms a part of '\textit{Ved\=anga}', a part essential to understand the \textit{Veda}s.

      	\textit{Chandash\=astra} by '\textit{Pi\.nga\d{l}a}' is the earliest treatise found on the Vedic Sanskrit meters. \textit{Pi\.nga\d{l}a} defines different meters on the basis of a sequence of what are called laghu and guru (short and long)  syllables and their count in the verse. The description and analysis of sequence of the laghu and guru syllables in a given verse is the major topic of \textit{Pi\.nga\d{l}a}'s work. He has described different sequences that can be constructed with a given number of syllables and has also named them. At the end of his book on '\textit{Chandash\=astra}', \textit{Pi\.nga\d{l}a}\superscript{3} gives rules to list all possible combinations of laghu and guru (L and G) in a verse with 'n' syllables, rules to find out the laghu-guru combinations corresponding to a given index, total number of possible combinations of 'n' L-G syllables and so on. In short \textit{Pi\.nga\d{l}a} describes the 'combinatorial mathematics' of meters in \textit{Chandash\=astra}. 
	Later around 8\superscript{th} century AD Shri \textit{Ked\=ar} \textit{Bhatt}\superscript{4} wrote '\textit{Vruttaratn\=akara}' a work on non-vedic meters. This seems to be independent of \textit{Pi\.nga\d{l}a}'s work, in the sense that it is not a commentary on \textit{Pi\.nga\d{l}a}'s work, and the last chapter gives the rules related to combinatorial mathematics which are totally different from \textit{Pi\.nga\d{l}a}'s approach. In the 13\superscript{th} century, \textit{Hal\=ayudha} in his '\textit{Mruta sanjivani}'\superscript{3} commentary on \textit{Pi\.nga\d{l}a}'s work, has again described the \textit{Pi\.nga\d{l}a}'s rules in great detail.

     	\textit{Pi\.nga\d{l}a}'s \textit{Chandash\=astra} contains 8 chapters. The 8\superscript{th}  chapter has 35 \textit{s\=utra}s of which the last 16 \textit{s\=utra}s from 8.20 to 8.35 deal with the algorithms related to combinatorial mathematics. \textit{Ked\=ar} \textit{Bhatt}'s '\textit{Vruttaratn\=akara}' contains 6 chapters, of which the 6\superscript{th}  chapter is completely devoted to algorithms related to combinatorial mathematics.

	Few words on the \textit{s\=utra} style of \textit{Pi\.nga\d{l}a} are in order. The \textit{s\=utra} style was prevalent during \textit{Pi\.nga\d{l}a}'s period. \textit{Asht\=adhy\=ayi} of \textit{P\=a\d{n}ini} is the classic example of \textit{s\=utra} style. A \textit{s\=utra} is defined as
\begin{center}
\textit{
	alp\=aksharam asandigdham s\=aravat vishvatomukham {\textbar}\\
	astobham anavadyam ca \textit{s\=utra}m \textit{s\=utra}vido viduH {\textbar}{\textbar}}
\end{center}

	A \textit{s\=utra} should contain minimum number of words (\textit{alp\=aksharam}), it should be unambiguous (\textit{asamdigdham}), it should contain essence of the topic which the \textit{s\=utra} is meant for (\textit{s\=aravat}), it should be general or should have universal validity(\textit{vishvatomukham}), it should not have any unmeaningful words (\textit{astobham}) and finally it should be devoid of any fault (\textit{anavadyam}).

	\textit{S\=utra}s are like mathematical formulae which carried a bundle of information in few words. They were very easy to memorise. They present a unique way to communicate algorithms or procedures verbally. The \textit{s\=utra} style was adopted by Indians in almost every branch of knowledge. For example, \textit{Pi\.nga\d{l}a}'s \textit{s\=utra}s are for combinatorics whereas Pa\d{n}ini's \textit{s\=utra}s are for language analysis. Another important feature of \textit{s\=utra} style is use of '\textit{anuvrutti}'. Generally all the \textit{s\=utra}s that deal with a particular aspect are clubbed together. To avoid any duplication utmost care had been taken to factor out the common words and place them at the appropriate starting \textit{s\=utra}.

Thus for example, if the following are the expanded \textit{s\=utra}s
\begin{flushleft}
		w\subscript{1} w\subscript{2} w\subscript{3} w\subscript{4}\\
		w\subscript{1} w\subscript{5} w\subscript{6}\\
		w\subscript{5} w\subscript{7} w\subscript{8}\\
		w\subscript{5} w\subscript{9}\\
		w\subscript{9} w\subscript{10} w\subscript{11} w\subscript{12}
\end{flushleft}
then the \textit{s\=utra} composer would put them as
\begin{flushleft}
		w\subscript{1} w\subscript{2} w\subscript{3} w\subscript{4}\\
		w\subscript{5} w\subscript{6}\\
		w\subscript{7} w\subscript{8}\\
		w\subscript{9}\\
		w\subscript{10} w\subscript{11} w\subscript{12}\\
\end{flushleft}
factoring out the words that are repeated in the following \textit{s\=utra}s.

One would then reconstruct the original forms by borrowing the words from the earlier \textit{s\=utra}s.  The context and the expectations provide the clues for borrowing. This process of borrowing or repeating the words from earlier \textit{s\=utra}s is known as '\textit{anuvrutti}'.  \textit{Pi\.nga\d{l}a} has used the \textit{s\=utra} style and also used '\textit{anuvrutti}'. \textit{Ked\=ar} \textit{Bhatt}'s \textit{Vruttaratn\=akara} contains \textit{s\=utra}s which are more verbose than that of \textit{Pi\.nga\d{l}a}'s, and does not use \textit{anuvrutti}.

In what follows, we take up each of the \textit{s\=utra}s from \textit{Pi\.nga\d{l}a}'s \textit{chandash\=astra} and explain its meaning and express it in modern mathematical language. We also examine the corresponding \textit{s\=utra} from \textit{Ked\=ar} \textit{Bhatt}'s \textit{Vruttaratn\=akara}, and compare the two algorithms.

\section{Algorithms}
	The algorithms that are described in \textit{Pi\.nga\d{l}a}'s  and \textit{Ked\=ar} \textit{Bhatt}'s work are

\begin{itemize}
\item  \textit{Prast\=araH}:   To get  all possible combinations (matrix) of n binary digits,
\item \textit{Nash\d{t}am}:  To recover the lost/missing  row in the matrix which is equivalent of getting a binary equivalent of a number,
\item \textit{Uddish\d{t}am}:  To get the row index of a given row in the matrix that is same as getting the value of a binary number,
\item \textit{Eka-dvi-\=adi-l-g-kriy\=a}: To compute  \superscript{n}C\subscript{r}, n being the number of syllables and r the number of  laghus (or \textit{gurus}), 
\item \textit{Samkhy\=a}:  To get the total number of n bit combinations; equivalent to computation of   2\superscript{n},
\item \textit{Adhva-yoga}: To compute the total combinations of \textit{chanda} (meters) ranging from 1 syllable to n syllables that is equivalent to computation of \begin{math}\displaystyle\sum_{i=1}^n 2^i\end{math}.
\end{itemize}

In addition to these algorithms, later commentators discuss an algorithm to get the positions of the r laghus in the matrix showing all possible combinations of n laghu-gurus. The corresponding structure is known as '\textit{pat\=ak\=a} \textit{prast\=ara}'.

\subsection{\textit{Prast\=araH}}
      We shall first give the \textit{Pi\.nga\d{l}a}'s algorithm followed by the \textit{Ked\=ar} \textit{Bhatt}'s.
\subsubsection{\textit{Pi\.nga\d{l}a}'s algorithm for \textit{Prast\=araH}}
	\textit{Prast\=araH} literally means expansion, spreading etc. From what follows it will be clear that by '\textit{prast\=ara}H', \textit{Pi\.nga\d{l}a} is talking about the matrix showing all possible combinations of n laghu-gurus. We know that there are 2\superscript{n} possible combinations of n digit binary numbers. So when we write all possible combinations, it will result into a   2\superscript{n}  *  n  matrix.

	The \textit{s\=utra}s in \textit{Pi\.nga\d{l}a}'s \textit{Chandash\=astra} are as follows:
\begin{table}[h]
\begin{center}
\begin{tabular}{|lr|} \hline
	dvikau glau 	&	8.20 \\ 
	misrau ca 	&	8.21 \\ 
	pruthagl\=a mishr\=aH & 	8.22 \\
	vasavastrik\=aH &	      	8.23 \\  \hline
\end{tabular}
\end{center}
\end{table}

\subsubsection{Explanation}

\begin{enumerate}
\item \begin{alltt}\textbf{dvikau glau}			8.20\end{alltt}
	This \textit{s\=utra} means \textit{prast\=ara} of 'one syllable(akshara)' has 2 possible elements viz. 'G  or  L'. So the 2\superscript{1}  * 1 matrix is 

\begin{math}
	\left[ \begin{array}{c} G\\ L    \end{array} \right] \end{math}

In  the boolean( 0 -1)  notation , if we put  0 for G  and  1 for L, we get
       
\begin{math} \left[ \begin{array}{c} 0\\ 1    \end{array} \right] \end{math}

\item \begin{alltt}\textbf{misrau ca}		         	8.21\end{alltt}
	“To get the \textit{prast\=ara} of two syllables, mix the above 1 syllable \textit{prast\=ara}   with itself”.  So we get 'G-L' mixed with 'G' followed by 'G-L' mixed with 'L'.
       
\begin{table}[h]
\caption{Mixed with G}
\begin{center}
\begin{tabular}{|l|l||l|} \hline
   G   &      	G                      &                   0  0 \\ 
   L   &	G                      &                   1  0 \\ \hline
\end{tabular}
\end{center}
\end{table}

\begin{table}[h]
\caption{mixed with L}
\begin{center}
\begin{tabular}{|l|l||l|} \hline
   G   &        L                      &                         0  1 \\ 
   L   &	L                      &                         1  1 \\ \hline
\end{tabular}
\end{center}
\end{table}

This will result in  the 2\superscript{2} * 2 matrix shown in table 3.
\begin{table}[h]
\caption{2 syllable combinations}
\begin{center}
\begin{tabular}{|lr||lr|} \hline
            G   &     G          &                  0 & 0 \\ 
            L   &     G          &                  1 & 0 \\ 
            G   &     L          &                  0 & 1 \\ 
            L	&     L          &                  1 & 1 \\  \hline
\end{tabular}
\end{center}
\end{table}

But in boolean notation we represent  all possible  2-digit numbers, in ascending order of their magnitude,  as in table 4.

\begin{table}[h]
\caption{2 digit Binary numbers}
\begin{center}
\begin{tabular}{|lr|} \hline
0 & 0  \\ 
0 & 1  \\ 
1 & 0  \\ 
1 & 1  \\  \hline
\end{tabular}
\end{center}
\end{table}

The Table (4)  is obtained by elementary transformation of exchange of columns from Table (3), or simply put it is just the mirror image of Table (4).
Now the question is, why the ancient Indian notation is as in  Table (3)  and not as in  (4).

This may be because of the practice of writing from LEFT to RIGHT. The characters uttered first are written to the left of those which are uttered later.

\item \begin{alltt} \textbf{pruthagl\=a mishr\=aH}		8.22 \end{alltt}
“To get the expansion of 3 binary numbers, again mix the G and L separately with the \textit{prast\=ara} of 2-syllables”. So we get a  2\superscript{3} * 3 matrix as in table 5.
                                                       
\begin{table}[h]
\caption{3 syllable \textit{prast\=ara}}
\begin{center}
\begin{tabular}{|lll|lll|lll|} \hline
\multicolumn{3}{c}{G-L-representation}&
\multicolumn{3}{c}{GL-01-conversion}&
\multicolumn{3}{c}{boolean notation}\\ \hline
            G &	G &	G   &  0 & 0 & 0 & 0 & 0 & 0\\
	    L &	G &	G   &  1 & 0 & 0 & 0 & 0 & 1\\
	    G &	L &	G   &  0 & 1 & 0 & 0 & 1 & 0\\
	    L &	L &	G   &  1 & 1 & 0 & 0 & 1 & 1\\  \hline
	    G &	G &	L   &  0 & 0 & 1 & 1 & 0 & 0\\
	    L &	G &	L   &  1 & 0 & 1 & 1 & 0 & 1\\
	    G &	L &	L   &  0 & 1 & 1 & 1 & 1 & 0\\
	    L &	L &	L   &  1 & 1 & 1 & 1 & 1 & 1\\ \hline
\end{tabular}
\end{center}
\end{table}

Again note that the modern (boolean notation) and ancient Indian notations (GL-01-conversion) are mirror images of each other.

\item vasavastrik\=aH		            8.23
\end{enumerate}
	This \textit{s\=utra} simply states that there are 8 (vasavaH) 3s (trik\=aH). 

Thus the first  of the 4 rules gives the terminating (or initial) condition. Second  rule tells how to generate a matrix for 2 bits from that of 1 bit. The third  rule states how to generate combinations for 3 bits, given combinations for 2 bits. Fourth rule describes the size of the matrix of 3 bits, and that's all. It is understood that this process (the 3\superscript{rd} rule)  is to be repeated again and again to get matrices of higher order.

\subsubsection{Recursive-ness of \textit{prast\=ara}}
To make it clear, we represent the matrix in step 1 as \\
\begin{math}
 A^{1}_{2*1} = \left[ \begin{array}{c} 0 \\ 1 \end{array}\right].\end{math} \\
Then the matrix in step 2 is\\
\begin{math}
 A^{2}_{4*2} = \left[ \begin{array}{cc}
0&0\\1&0\\0&1\\1&1 \end{array}\right] 
=\left[ \begin{array}{cc} A^{1}_{2*1} &0_{2*1}\\
A^{1}_{2*1} & 1_{2*1} \end{array} \right].  \end{math} 
\\
where \begin{math}O_{m*n}\end{math} is a matrix with all elements equal to 0 and \begin{math}1_{m*n}\end{math} is a matrix with all elements equal to 1.\\

Continuing further, the matrix in step 3 is \\
\begin{math}
A^{2}_{4*2} = \left[ \begin{array}{ccc} 0&0&0 \\ 1&0&0 \\0&1&0 \\ 1&1&0 \\ 0&0&1\\ 1&0&1 \\ 0&1&1 \\ 1&1&1 \end{array} \right] = \left[  \begin{array}{cc}A^{2}_{4*2} 0_{4*1}\\
   A^{2}_{4*2} 1_{4*1} \end{array} \right]. \end{math}\\

The generalisation of this leads to \\
\begin{math}
A^{n}_{2^{n}*n} = \left[ \begin{array}{cc} A^{n-1}_{2^{n-1}*(n-1)} 0_{2^{n-1}*1}\\
				   A^{n-1}_{2^{n-1}*(n-1)} 1_{2^{n-1}*1} \end{array} \right]\end{math}

	We notice that the algorithm for generation of all possible combinations of n bit binary numbers is thus 'recursive'.
	
\subsubsection{\textit{Ked\=ar} \textit{Bhatt}'s algorithm for Prast\=araH}
	\textit{Ked\=ar} \textit{Bhatt} in his \textit{Vruttaratn\=akara} has given another algorithm to get the \textit{prast\=ara} for a given number of bits. His algorithm goes like this:

\begin{alltt}
\begin{center}
	p\=ade sarvagurau \=ady\=at laghu nasya guroH adhaH \begin{math}|\end{math}
	yath\=a-upari tath\=a sesham bh\=uyaH kury\=at amum vidhim \begin{math}||\end{math}
	\=une dady\=at gur\=un eva y\=avat sarve laghuH bhavet \begin{math}|\end{math}
	\textit{prast\=ara}H ayaM sam\=akhy\=ataH chandoviciti vedibhiH \begin{math}||\end{math}

\end{center}
\end{alltt}

"In the beginning all are gurus(G)(\textit{p\=ade sarvagurau}) . In the second line, place a laghu(L) below the first G of the previous line(\textit{\=ady\=at laghu nasya guroH adhaH}) . Copy the remaining  right  bits as in the above line(\textit{yath\=a-upari tath\=a sesham}) . Place Gs in all the remaining places to the left (if any) of the 1\superscript{st}  G bit(\textit{\=une dady\=at gur\=un eva}). Repeat this till all of them become laghu(\textit{y\=avat sarve laghuH bhavet}) . This process is known as '\textit{prast\=ara}'."

Here is an example, explaining the above algorithm:

\begin{table}[h]
\caption{Kedar Bhatt's algorithm}
\begin{center}
\begin{tabular}{|lll|p{6cm}|} \hline
           G&	G&	G & start with all Gs\\
	\textbf{L}&	G&	G & place L below 1\superscript{st}  G  of above line, copy remaining  right bits viz. G G as in the above line\\
	G&	\textbf{L}&	G & place L below 1\superscript{st}  G of above line, copy remaining  right bit viz. G as in the above line,  and G in the remaining place to the left of L\\
	\textbf{L}&	L&	G & place L below 1\superscript{st}  G of above line, copy remaining  right bits viz. L G as in the above line\\
	G&	G&	\textbf{L} & place L below 1\superscript{st} G of above line and G in  remaining places  to the left\\
	\textbf{L}&	G&	L & place L below 1\superscript{st}  G of above line, copy remaining  right bits  viz. G L as in the above line\\
	G&	\textbf{L}&	L &place L below  1\superscript{st}  G of above line, copy remaining right bit viz. L as in above line and G in the remaining place to  the left\\
	\textbf{L}&	L&	L &place L below  1\superscript{st}  G of above line, copy remaining remaining right bits viz. L L  as in the above line and stop the process since all are Ls.\\ \hline
\end{tabular}
\end{center}
\end{table}
\pagebreak
	If we compare \textit{Pi\.nga\d{l}a}'s method with that of \textit{Ked\=ar} \textit{Bhatt}'s , we note that the first one is a recursive, whereas the second one is an iterative one.
	Compare this with the well-known definitions of factorials in modern notation. We can define factorial  in two different ways:
	
\begin{center}
\begin{verbatim}
	n! = n * (n-1)!
	1! = 1

	OR
	
	n! = 1*2*3*...*n .
\end{verbatim}
\end{center}

Thus if one uses the first definition to get a factorial of say 4, one needs to know how to get the factorial of 3; to get the factorial of 3, one in turn should know how to get factorial of 2, etc.\\
Similarly according to \textit{Pi\.nga\d{l}a}'s algorithm, to write a \textit{prast\=ara} for 4 syllables, one needs to write a \textit{prast\=ara} for 3 syllables, and to do so in turn one should write a \textit{prast\=ara} for 2 syllables, and so on.\\
On the other hand, using the second definition, one can get 4! just by multiplying 1,2,3 and 4. One need not go through the whole process of finding other factorials. Similarly \textit{Ked\=ar Bhatt} describes an algorithm where one can write the \textit{prast\=ara} for say 4 syllables directly without knowing what the \textit{prast\=ara} for 3 syllables is.\\
The  algorithm to obtain \textit{prast\=ara}, as given by \textit{Pi\.nga\d{l}a}, is similar to the recursive definition of factorial, whereas the one  given by \textit{Ked\=ar} \textit{Bhatt}, is similar to the iterative definition of factorial.

\subsection{\textit{Nash\d{t}am}}
	In ancient days, the \textit{prast\=ara} (or matrix) used to be written on the sand, and hence there was possibility of getting a row erased. The next couple of \textit{s\=utra}s (8.24 and 8.25) are to recover the lost (or disappeared or vanished) row (nash\d{t}a) from the matrix. If one knows \textit{Ked\=ar} \textit{Bhatt}'s algorithm then the lost  row can be recovered easily from the previous or the next one. But \textit{Pi\.nga\d{l}a} did not have an iterative description and hence he has given a separate algorithm to find the 'lost' row. In different words, getting a 'lost' row is conceptually equivalent to getting guru-laghu combination (i.e. the binary equivalent) of  the row index.

	The \textit{s\=utra}s are as follows:

\begin{table}[h]
\begin{center}
\begin{tabular}{|lr|} \hline
	l-arddhe      &		(8.24) \\ 
	sa-eke-ga      &		(8.25)  \\ \hline
\end{tabular}
\end{center}
\end{table}

"In case the given number can be halved (without any remainder), then write 'L', else add one and then halve it and write 'G'". For example, suppose we want to get the 'laghu-guru' combination in the 5\superscript{th} row of the 3-akshara matrix. We start with the given row-index i.e. 5. Since it is an odd number, add 1 to it and write 'G'. After dividing 5+1(=6) by 2, we get 3. Again this is an odd number, and hence we add 1 to it, and write 'G'. After dividing 3+1(=4) by 2 we get 2. Since it is an even number we write 'L'. Once we get desired number of bits (in this case 3), the process ends:

\begin{table}[h]
\caption{\textit{Nas\d{t}am}}
\begin{center}
\begin{tabular}{|lrr|l|} \hline
	5 &-\begin{math}>\end{math}& (5+1)/2=3&		G \\
	3 &-\begin{math}>\end{math}& (3+1)/2=2&		G G \\
	2 &-\begin{math}>\end{math}& 2/2=1	&	G G L  \\ \hline
\end{tabular}
\end{center}
\end{table}
So the 5\superscript{th} row in the \textit{prast\=ara} (matrix) of 3 bits is G G L.  The algorithm may be written as a recursive function as follows:

\begin{verbatim}
	Get_Binary(n) =
		Print L ; Get_Binary(n / 2), if n is even,
		Print G ; Get_Binary(n+1 / 2), if n is odd,
		Print G; if n=1. (terminating condition)
\end{verbatim}
Thus this algorithm gives a method to convert a binary equivalent of a given number.

\subsubsection{Difference between \textit{Pi\.nga\d{l}a}'s method and Boolean method}
Let us compare this conversion with the modern method. The boolean method is illutrated in table 8.

\begin{table}[h]
\caption{Boolean conversion to binary}
\begin{center}
\begin{tabular}{|llr|} \hline
	5		& remainder & \\
	5 / 2 = 2 & 1 &	\^{} \\
	2 / 2 = 1	&	0& 	\begin{math}|\end{math}\\
	1 / 2 = 0	&	1& ----\begin{math}>|\end{math} \\ \hline
\end{tabular}
\end{center}
\end{table}

Hence the binary equivalent of 5 is 101. If we replace G by 0 and L by 1 in 'G G L' we get  0 0 1. We have seen earlier that the numbers in modern and Indian method are mirror images, so after taking the mirror image of '0 0 1' we get '1 0 0'. Thus, by \textit{Pi\.nga\d{l}a}'s method we get the equivalent of 5 as '1 0 0' whereas by modern method, we get 5=101\subscript{2}. Why is this difference? This difference is attributed to the fact that the counting in \textit{Pi\.nga\d{l}a}'s method starts with '1'. In other words, 1 is represented as '0 0 0' in \textit{Pi\.nga\d{l}a}'s method, and not as '0 0 1'.

	Thus we notice two major differences between the \textit{Pi\.nga\d{l}a}'s method and the modern representation of binary numbers viz.\\
in \textit{Pi\.nga\d{l}a}'s system,
\begin{itemize}
\item as has been initially observed by Nooten\superscript{2}, the numbers are written with the higher place value digits to the right of lower place value digits, and
\item the counting starts with 1.
\end{itemize} 
\subsection{\textit{Uddish\d{t}am}}
	The third algorithm is to obtain position of the desired (\textit{uddishta}) row in a given matrix, without counting its position from the top, i.e. to get the row index corresponding to a given combination of G and Ls. Thus this is the inverse operation of \textit{nash\d{t}am}.
       	Both \textit{Pi\.nga\d{l}a} as well as \textit{Ked\=ar} \textit{Bhatt} have given algorithms for \textit{uddish\d{t}am}.

\subsubsection{\textit{Pi\.nga\d{l}a}'s algorithm for \textit{uddish\d{t}am}}
	Two \textit{s\=utra}s viz. (8.26) and (8.27) from \textit{Pi\.nga\d{l}a}'s \textit{Chandash\=astra} describe this algorithm.
The \textit{s\=utra}s are as follows:

\begin{table}[h]
\begin{center}
\begin{tabular}{|lr|} \hline
	pratilomaguNam dviH-l-\=adyam	&	(8.26) \\
	tat\=aH-gi-ekaM jahy\=at	&		(8.27) \\ \hline
\end{tabular}
\end{center}
\end{table}

We first see the meaning of these \textit{s\=utra}s followed by an example.
	The first \textit{s\=utra} states that in the reverse order(\textit{pratilomaguNam}), starting from the 1\superscript{st} laghu(\textit{l-\=adyam}), multiply by 2(\textit{dviH}). The second \textit{s\=utra} states that if the syllable is guru(\textit{gi}), subtract one (\textit{ekaM jahy\=at}) (after multiplying by 2). Here we also note the use of '\textit{anuvrutti}'. The word \textit{dviH} is not repeated in the following \textit{s\=utra}, but should be borrowed from the previous \textit{s\=utra} .
Since it is not mentioned what the starting number should be, we start with 1.

We illustrate this with an example. Let the input sequence be 'G L G'. Table 9 describes the application of the above sutras.

\begin{table}[h]
\caption{\textit{uddis\d{t}am}}
\begin{center}
\begin{tabular}{|lll|p{5cm}|} \hline
	G &	L &	G & remark\\ \hline
	 &	1 &	 &	(start with 1\superscript{st} L from the right, starting number 1)\\
 &		2 & &		(multiply by 2)\\
	2 & & &			(continue with the previous result i.e. 2)\\
	4   &&&			(multiply by 2)\\
	3   &&&			(subtract 1, since it is guru) .\\ \hline
\end{tabular}
\end{center}
\end{table}

Thus the row 'G L G' is in the 3\superscript{rd} position in the \textit{prast\=ara} of 3 bits.  It is clear that this set of rules thus gives the row index of a row in the \textit{prast\=ara} matrix.

The algorithm may be written formally as in table 10. \\
\begin{table}[h]
\caption{algorithm for Base 2}
\begin{center}
\begin{tabular}{|lp{7cm}|} \hline
S\subscript{i} & = 1 where 1\superscript{st} laghu occurs in the i+1\superscript{th} position from right.\\
S\subscript{i+1}&  = 2 * S\subscript{i}     if i+1\superscript{th} position has L,\\
                &  = 2 * S\subscript{i} - 1 if i+1\superscript{th} position has G,\\
 & where S\subscript{i} denotes index of i\superscript{th} digits from the right.\\ \hline
\end{tabular}
\end{center}
\end{table}

This set of rules further can be extended to get the decimal value of a number in any base B as shown in table 11.

\begin{table}[h]
\caption{algorithm:for Base B}
\begin{center}
\begin{tabular}{|lp{7cm}|} \hline
S\subscript{i} & = 1 where 1\superscript{st} non-digit zero occurs in the i+1\superscript{th} position. (The counting for i starts with 1, and goes from the right digit with highest place value to the lowest place value)\\
S\subscript{i+1}&  = B * S\subscript{i}     if i+1\superscript{th} position has B-1,\\
                &  = B * S\subscript{i} - D\superscript{'}\subscript{i+1}, otherwise,\\
 & where D\superscript{'}\subscript{i+1} stands for the B-1\superscript{'}s complements of i+1\superscript{th} digit, \\&and S\subscript{i} denotes index of i\superscript{th} digit from the right.\\ \hline
\end{tabular}
\end{center}
\end{table}

According to this algorithm, value of the decimal number 789 can be calculated as shown in the table 12.
\begin{table}[h]
\caption{Example: base 10}
\begin{center}
\begin{tabular}{|lp{5cm}|} \hline
	S0  &=  1\\
	S1  &=  10 * 1 - 2  ( 9's complement of 7)\\
      	    &  =  8\\
	S2 &=   8 * 10 - 1 (9's complement of 8)\\
           &      = 79\\
	S3 &= 79 * 10 \\
           &      = 790 .\\ \hline
\end{tabular}
\end{center}
\end{table}
Thus the position of 789 in the decimal place value system (where 0 is the 1\superscript{st} number 1 is the 2\superscript{nd} and so on) is 790.
Further note that S1(=8) gives the value of single digit 7, S2(=79) gives the index of 2 digits 78.

\subsubsection{\textit{Ked\=ar} \textit{Bhatt}'s algorithm for \textit{uddish\d{t}am}}
	The \textit{Ked\=ar} \textit{Bhatt}'s version of '\textit{uddish\d{t}am}' differs from that of \textit{Pi\.nga\d{l}a}. \textit{Ked\=ar} \textit{Bhatt}'s version goes like this:

\begin{center}
\begin{alltt}
	uddish\d{t}am dvigun\=an \=ady\=an upari ank\=an sam\=alikhet \begin{math}|\end{math}
	laghusth\=a ye tu tatra ank\=ataiH sa-ekaiH mishritaiH bhavet \begin{math}|\end{math}\end{alltt}
\end{center}

"To get the row number corresponding to the given laghu guru combination, starting from the first, write double (the previous one) on the top of each laghu-guru. Then all the numbers on top of laghu are added with 1. (Since the starting number is not mentioned, by default, we start with 1)."

We illustrate this with an example.

Let the row be 'G L L'.

We start with 1, write it on top of 1\superscript{st} G. Then multiply it by 2, and write 2 on top of L, and similarly (2*2=) 4 on top of the last L.

\begin{table}[h]
\caption{\textit{uddis\.tam}}
\begin{center}
\begin{tabular}{|lll|} \hline
	1  &	2  &	4 \\
	G  &	L  &	L \\ \hline

\end{tabular}
\end{center}
\end{table}
Then we add all the numbers which are on top of L viz. 2 + 4. To this we add 1. So the row index of the given L-G combination is 7.

	In boolean mathematical notation, 'G L L' stands for 1 1 0 (after mirror image). This is equivalent to decimal 6. The difference of 1 is attributed to the fact that row index is counted from 1, as pointed out earlier.

\subsection{eka-dvi-\=adi l-g kriy\=a}
       \textit{Ked\=ar} \textit{Bhatt}'s work describes explicit rules to get the number of combinations of 1L, 2L, etc. (1G, 2G, etc.) among all possible combinations of n L-Gs. In other words, it gives a procedure to calculate \superscript{n}C\subscript{r} . \textit{Pi\.nga\d{l}a}'s \textit{s\=utra} is very cryptic and it is only through \textit{Hal\=ayudha}'s commentary on it, one can interpret the \textit{s\=utra}as a 'meru' which resembles the Pascal's triangle. we first give \textit{Ked\=ar} Bhatta's algorithm followed by \textit{Pin\.nga\d{l}a}'s.

\subsubsection{\textit{Ked\=ar} \textit{Bhatt}'s algorithm}
	The procedure for \textit{eka-dvi-adi-l-g-kriy\=a} in \textit{vruttaratn\=akara} is described as follows:

\begin{center}
	\begin{alltt}
	varn\=an vrutta bhav\=an sa-ek\=an auttar\=ardhayataH sthit\=an \textbar
	ek\=adikramataH ca et\=an upari-upari nikshipet \textbar\textbar
	up\=antyato nivartet tyajat na ekaikam \=urdhavataH \textbar
	upari \=ady\=at guroH evaM eka-dvi-\=adi-l-g-kriy\=a \textbar\textbar 
\end{alltt}
\end{center}

Whatever the given number of syllables is, write those many 1s from the left to right as well as top to bottom. Then in the 1\superscript{st} row, add the number in the top(previous) row to its left occupant, and continue this process leaving the last number. Continue this process for the remaining rows. The last number in the 1\superscript{st} column stands for all gurus. The last number in the second column  stands for one laghu, the one in the third column for two laghus, and so on.

We explain this algorithm by an example.

Let the number be 6.
Write 6 1s horizontally as well as vertically as below. Elements are populated rowwise by writing the sum of numbers in immediately preceeding row and column.
\begin{table}[h]
\caption{\textit{Meru} aka Pascal's Triangle}
\begin{center}
\begin{tabular}{|lllllll|} \hline
    &	1   &	1    &	1  &	1  &	1   &	\textbf{1} \\
1   &	2   &	3    &	4  &	5  &	\textbf{6} & \\
1   &	3   &	6    &	10  &	\textbf{15} & & \\
1   &	4   &	10   &	\textbf{20}  &&	 & \\
1   &	5   &	\textbf{15}     &&& & \\
1   &	\textbf{6}   &  &&& & \\
\textbf{1}   & &  &&& & \\ \hline
\end{tabular}
\end{center}
\end{table}

The numbers 1, 6, 15, 20, 15, 6, 1 give number of combinations with all gurus, one laghu, two laghus, three laghus, four laghus, five laghus, and finally all laghus.

	We see the striking similarity of this expansion with the Pascal's triangle. This process describes the method of getting the number of combinations of r from n viz. \superscript{n}C\subscript{r} .  This triangle is termed as \textit{meru} (literal: hill) in the Indian literature.

\subsubsection{\textit{Pi\.nga\d{l}a}'s algorithm}

{\textit{Pi\.nga\d{l}a}'s \textit{s\=utra}s are\\
\begin{table}[h]
\begin{center}
\begin{tabular}{|lr|} \hline
	\textit{pare p\=ur\d{n}am}      &		(8.34) \\ 
	\textit{pare p\=ur\d{n}am iti}      &		(8.35)  \\ \hline
\end{tabular}
\end{center}
\end{table}

The sutra 8.34 literally means, "complete it using the two far ends \textit{pare}". Only from the \textit{Hal\=ayudha}'s commentary it becomes clear that this \textit{s\=utra} means: Start with '1' in a cell. Below this cell draw two cells, and so on. Then fill all the cells which are at the far ends, in each row, by 1s. This results in figure 15.
\\
\begin{table}[h]
\begin{center}
\caption{\textit{Meru} construction step-1}
\begin{tabular}{|lllllllllllllll|} \hline
    &&&&&&&1&&&&&&& \\
    &&&&&&1&&1&&&&&& \\
    &&&&&1&&&&1&&&&& \\
    &&&&1&&&&&&1&&&& \\
    &&&1&&&&&&&&1&&& \\
    &&1&&&&&&&&&&1&& \\
    &1&&&&&&&&&&&&1& \\
    1&&&&&&&&&&&&&&1 \\ \hline
\end{tabular}
\end{center}
\end{table}

Next \textit{s\=utra} says, complete a cell using above two cells, again filling the far end cells. Thus resulting in table 16. Repeating this process we get table 17, and we see that the repeatition leads to the building of meru, or pascal's triangle. \\

\begin{table}[h]
\caption{\textit{Meru} construction step-2}
\begin{center}
\begin{tabular}{|lllllllllllllll|} \hline
    &&&&&&&1&&&&&&& \\
    &&&&&&1&&1&&&&&& \\
    &&&&&1&&2&&1&&&&& \\
    &&&&1&&3&&3&&1&&&& \\
    &&&1&&4&&&&4&&1&&& \\
    &&1&&5&&&&&&5&&1&& \\
    &1&&6&&&&&&&&6&&1& \\
    1&&7&&&&&&&&&&7&&1 \\ \hline
\end{tabular}
\end{center}
\end{table}

\begin{table}[h]
\caption{\textit{Meru} construction step-3}
\begin{center}
\begin{tabular}{|lllllllllllllll|} \hline
    &&&&&&&1&&&&&&& \\
    &&&&&&1&&1&&&&&& \\
    &&&&&1&&2&&1&&&&& \\
    &&&&1&&3&&3&&1&&&& \\
    &&&1&&4&&6&&4&&1&&& \\
    &&1&&5&&10&&10&&5&&1&& \\
    &1&&6&&15&&&&15&&6&&1& \\
    1&&7&&21&&&&&&21&&7&&1 \\ \hline
\end{tabular}
\end{center}
\end{table}

\subsubsection{Bhaskar\=acharya's method of obtaining \textit{meru}}
	There are other ways of obtaining this 'meru' described in Indian literature. For example Bhaskar\=acharya-II -{} the 12\superscript{th} century Indian mathematician - in his L\=il\=avati\superscript{5} gives following procedure for obtaining n\superscript{th} row of the meru. Write the numbers from 1 till n in an order(row A in table 18). Write the same numbers above this row, but this time in reverse order(row B in table 18). Now start with 1, and proceed as in the following algorithm:

\begin{table}
\caption{Binary Coefficients}
\begin{center}
\begin{tabular}{|llllllll|} \hline
i&	6&	5&	4&	3&	2&	1&	0\\
C&	1&	6&	15&	20&	15&	6&	1\\
B&       &	6&	5&	4&	3&	2&	1\\
A&       &	1&	2&	3&	4&	5&	6\\ \hline
\end{tabular}
\end{center}
\end{table}

\begin{center}
	C\subscript{0} = 1 \\
	A\subscript{i} = n - i \\
	B\subscript{i} = i +1 \\
	C\subscript{i+1} = C\subscript{i}  * A\subscript{i} / B\subscript{i} . \\
                    
\end{center}
Or, in other words,

\begin{center}
	\superscript{n}C\subscript{0} = 1\\
	\superscript{n}C\subscript{r+1} = \superscript{n}C\subscript{r} * (n-r)/(r+1) .\\

\end{center}
This is another instance of recursive definition.

\subsection{\textit{Sankhy\=a}}
	'Sankhy\=a' stands for the number of possible combinations of n bits. \textit{Pi\.nga\d{l}a} and \textit{Ked\=ar} give an algorithm to calculate 2\superscript{n} , given n. The algorithms differ as in earlier cases. \textit{Ked\=ar} \textit{Bhatt} uses the results of previous operations (\textit{uddish\d{t}am} and \textit{eka-dvi-\=adi-l-g-kriya}), whereas \textit{Pi\.nga\d{l}a} describes a totally independent algorithm.

\subsubsection{\textit{Ked\=ar} \textit{Bhatt}'s algorithm for finding the value of 'Sankhy\=a'}
	\textit{Ked\=ar} \textit{Bhatt} gives the following \textit{s\=utra} in his 6\superscript{th} chapter of the book '\textit{vruttaratn\=akara}'

\begin{alltt}
\begin{center}
	l-g-kriy\=anka sandohe bhavet sankhy\=a vimishrite \begin{math}|\end{math}
	uddish\d{t}a-anka sam\=ah\=ar\=aH sa eko v\=a janayedim\=am \begin{math}||\end{math}
\end{center}
\end{alltt}

This \textit{s\=utra} says, one can get the total combinations in two different ways:\\
a) by adding the numbers of \textit{eka-dvi-\=adi-l-g-kriy\=a} , or\\
b) by adding the numbers at the top in the \textit{uddish\d{t}a kriy\=a} and then adding 1 to it.

So for example, to get the possible combinations of 6 bits,\\
\begin{itemize}
\item the numbers in the \textit{eka-dvi-\=adi-l-g-kriy\=a} are 1,6,15,20,15,6,1 (see table 14). Adding these we get\\
	1 + 6 + 15 + 20 + 15 + 6 + 1 = 64.\\
Therefore, there are 64 combinations of 6 bits.
\item The \textit{uddish\d{t}a} numbers in case of 6 bits are\\
	1,2,4,8,16,32 \\
and adding all these and then 1 to it, we get\\
	1+2+4+8+16+32+1 = 64.
\end{itemize}
From this it is obvious that \textit{Ked\=ar} \textit{Bhatt} was aware of the following two well-known formulae.

	2\superscript{n} = \begin{math}\displaystyle\sum_{r=0}^{n}\end{math} \superscript{n}C\subscript{r}      (Sum of the numbers in \textit{eka-dv-\=adi-l-g-kriy\=a})

and 

	2\superscript{n} = \begin{math}\displaystyle\sum_{i=0}^{n-1} 2^{i}\end{math}  + 1  (sum of \textit{uddish\d{t}a} numbers +1) .
\subsubsection{\textit{Pi\.nga\d{l}a}'s algorithm for finding the value of '\textit{sankhy\=a}'}
	\textit{Pi\.nga\d{l}a}'s description goes like this:

\begin{table}[h]
\begin{center}
\begin{tabular}{|lr|} \hline
	dviH arddhe	&		(8.28) \\
	rūpe sh\=unyam	&		(8.29) \\
	dviH sh\=unye	&		(8.30) \\
	t\=avadardhe tadgunitam &	(8.31)  \\ \hline
\end{tabular}
\end{center}
\end{table}

If the number is divisible by 2\{\textit{arddhe}\}, divide by 2 and write 2\{\textit{dviH}\}.
If not, subtract 1\{\textit{r\=upe}\}, and write 0\{\textit{sh\=unyam}\}.
If the answer were 0\{\textit{sh\=unya}\}, multiply by 2\{\textit{dviH}\}, and if the answer were 2\{\textit{arddhe}\}, multiply \{\textit{tad gu\d{n}itam}\} by itself \{\textit{t\=avad}\}.

So for example, consider 8.
\begin{verbatim}
	8
	4	2  (if even, divide by 2 and write 2)
	2	2  (if even, divide by 2 and write 2)
	1	2  (if even, divide by 2 and write 2)
	0	0 (if odd, subtract 1 and write 0) .
\end{verbatim}
Now start with the 2\superscript{nd} column, from bottom to top.
\begin{verbatim}
		0	1*2 = 2          (if 0, multiply by 2)
		2	2^2 = 4          (if 2, multiply by itself)
		2	4^2 = 16        (if 2, multiply by itself)
		2	16^16 = 256   (if 2, multiply by itself).
\end{verbatim}

This algorithm may be expressed formally as
\begin{verbatim}

        power2(n)   =   [power2(n/2)] ^  2  if n is even,
                    =   power2(n-1/2) * 2, if n is odd,
                    =   1, if n = 0 .
		       
\end{verbatim}
Note that the results after each call of the function are 'stacked' and may also be treated as 'tokens' carrying the information for the next action (whether to multiply by 2 or to square). It still remains unclear to the author which part of the \textit{s\=utra} codes information about  'stack'. Or, in other words, how does one know that the operation is to be done in reverse order? There is no information about this in the \textit{s\=utra}s anywhere either explicit or implicit. This algorithm of calculating n\superscript{th} power of 2 is a recursive algorithm and its complexity is O(log\subscript{2}n), whereas the complexity of  calculating power by normal multiplication is O(n). Knuth\superscript{6} has referred to this algorithm as a 'binary method' (Knuth, pp 399). 

\subsection{adhvayoga}
\textit{Pi\.nga\d{l}a}'s \textit{s\=utra} is\\

\begin{table}[h]
\begin{center}
\begin{tabular}{|lr|} \hline
        dviH dviH \=unam tad ant\=an\=am  &  8.32 \\ \hline
\end{tabular}
\end{center}
\end{table}

This algorithm gives the sum (\textit{yoga}) of all the chandas (\textit{adhva}) with number of syllables less than or equal to n. The sutra literally means to get \textit{adhvayoga}, multiply the last one (\textit{tat ant\=an\=am}) by 2 (\textit{dviH}) and then subtract 2 (\textit{dviH Unam}). That is\\
\begin{center}
\begin{math}\displaystyle\sum_{i=1}^{n} 2^{i}\end{math} = 2\superscript{n} * 2 - 2 = 2 \superscript{n+1} - 2\\
\end{center}
\subsection{Finding the position of  all combinations of r guru (laghu) in a \textit{prast\=ara} of n bits}
	This is an interesting algorithm found only in commentaries on \textit{Ked\=ar} \textit{Bhatt}'s work\superscript{4}. We have not been able to trace the origin of this algorithm. This is to find the positions of combinations involving 1 laghu, 2 laghu, etc. in the n bit \textit{prast\=ara}. For example, in the 2 bit \textit{prast\=ara} shown in table (2), we see that there is only one combination with both Gs, and it occurs in the 1\superscript{st} position. There are 2 combinations of 1G (or 1 L), and they occur at the 2\superscript{nd} and the 3\superscript{rd} positions. Finally there is only one combination of 2 Ls, and it occurs at the 4\superscript{th} position. The following algorithm describes a way to get these positions without writing down the \textit{prast\=ara}.

\subsubsection{Algorithm to get positions of r laghu(guru) in a \textit{prast\=ara} of n laghu-gurus}
We will give an algorithm to populate the matrix A such that the j\superscript{th} column of A will have positions of the rows in \textit{prast\=ara} with j laghus. It follows that the total number of elements in j\superscript{th} column will be \superscript{n}C\subscript{j} .\\
\begin{enumerate}
\item Write down 1,2,4,8,...,2\superscript{n} in the 1\superscript{st} row.\\
    
       A[0,i] = 2\superscript{i} ,     0 \begin{math}\le\end{math} i \begin{math}\le\end{math} n .\\

\item The 2\superscript{nd} column of elements is obtained by the following operation:\\
    A[1,i] = A[0,0] + A[0,i] , 1 \begin{math}\le\end{math}  i \begin{math}\le\end{math}  n, and A[i,j] \begin{math}<\end{math} 2\superscript{n} .\\

\item The remaining columns (3\superscript{rd} onwards) are obtained as follows:\\
 For each of the elements A[k-1,j] in the k\superscript{th} column, do the following:    \\
\begin{alltt}
	A[k,m] = A[k-1,j] + A[0,i],  k \begin{math}\le\end{math} i \begin{math}\le\end{math}  n+1, \\
        if A[k,m] does not occur in the so-far-populated matrix, and\\
        A[k,m] \begin{math}<\end{math} 2\superscript{n},  and \\
        0 \begin{math}\le\end{math} j \begin{math}\le\end{math} \superscript{n}C\subscript{j}, and \\
        0  \begin{math}\le\end{math} m \begin{math}\le\end{math} \superscript{n}C\subscript{l} . 
\end{alltt}
\end{enumerate}

The 1\superscript{st} column gives positions of rows with all gurus;\\
2\superscript{nd} column gives position of rows with 1 laghu, and remaining gurus;\\
3\superscript{rd} column gives position of rows with 2 laghu and remaining gurus, and so on.\\
The last column gives the position of rows with all laghus.\\

Following example will illustrate the procedure.\\
Suppose we are interested in the positions of different combinations of laghus and gurus in the \textit{prast\=ara} of 5 bits. We start with the powers of 2 starting from 0 till 5 as the 1\superscript{st} row.\\

\begin{verbatim}
	1	2	4	8	16	32 .
\end{verbatim}

To get the 2\superscript{nd} column, we add 1 (A[0,0]) to the remaining elements in the 1\superscript{st} row (see table 19).
               
\begin{table}[h]
\caption{populating 2\superscript{nd} column}
\begin{center}
\begin{tabular}{|l|l|l|l|l|l|p{3cm}|} \hline
	1&	2&	4&	8&	16&	32&\\
	&	3&	&	&	&	&	(1 + 2)\\
	&	5&	&	&	&	&	(1 + 4)\\
	&	9&	&	&	&	&	(1 + 8)\\
	&	17&	&	&	&	&	(1 + 16) .\\ \hline
\end{tabular}
\end{center}
\end{table}

Thus, the 2\superscript{nd} column gives the positions of rows in the \textit{prast\=ara} of 5 bits with 1 laghu and remaining (4) gurus.

To get the 3\superscript{rd} column, we add 2 (A[1,0]) to the remaining elements(A[0,j]; j \begin{math}>\end{math} 1) of the 1\superscript{st} row (see table 20).

\begin{table}[h]
\caption{populating 3\superscript{rd} column}
\begin{center}
\begin{tabular}{|l|l|l|l|l|l|p{3cm}|} \hline
	1&	2	&4	&8	&16	&32 &\\
	&	3	&6	&	&	&&	(2 + 4)\\
	&	5	&10	&	&	&&	(2 + 8)\\
	&	9	&18	&	&	&&	(2 + 16)\\
	&        17 &&&&&.\\		\hline
\end{tabular}
\end{center}
\end{table}
We repeat this for other elements in the 2\superscript{nd} column (A[i,1]; i\begin{math} > \end{math}0) as in Table 21.

\begin{table}[h]
\caption{populating 3\superscript{rd} column contd}
\begin{center}
\begin{tabular}{|l|l|l|l|l|l|p{3cm}|} \hline
	1&	2 &	 4&	8&	16&	32& \\
	&	3 &	 6&	&	&	&	(2 + 4)\\
	&	5 &	10&	&	&	&	(2 + 8)\\
	&	9 &      18&	&	&	&	(2 + 16)\\
	&       17&  	 7&	&	&	&	(3 + 4)\\
		  &&      11&	&	&	&	(3 + 8)\\
	&&		19&	&	&	&	(3 + 16) \\
	&&		13&	&	&	&	(5 + 8) [ 5+4=9 already exists in the matrix, and hence ignored]\\
	&&		21&	&	&	&	(5+16)\\
	&&		25&	&	&	&	(9+16) [ 9+4, 9+8 are ignored] .\\ \hline
\end{tabular}
\end{center}
\end{table}
The 3\superscript{rd} column gives positions of rows with 2 laghus in the 5 bit expansion.
We repeat this procedure till all the columns are exhausted. The final matrix will be as in table 22.

\begin{table}[h]
\caption{\textit{pat\=ak\=a}}
\begin{center}
\begin{tabular}{|l|l|l|l|l|l|} \hline
	1 &	2 &	4 &	8 &	16 &	32  \\ 
	&	3 &	6 &	12 &	24 &  \\ 
	&	5 &	10 &	20 &	28 &  \\ 
	&	9 &	18 &	14 &	30 &  \\ 
	&	17 &	7 &	22 &	31 & 	 \\ 
& &			11 &	26 & & \\ 
& &			19 &	15 & & \\ 
& &			13&	23& &\\ 
& &			21&	27& &\\ 
& &			25&	29& & \\  \hline
\end{tabular}
\end{center}
\end{table}

Since this is in the form of '\textit{pat\=ak\=a}'(which literally means a flag\footnote{Indian flags used to be triangular in shape}), it is also called as '\textit{pat\=ak\=a} \textit{prast\=ara}'. 
Thus the Indian mathematicians have gone one step ahead of the modern mathematicians and not only gave the algorithms to find \superscript{n}C\subscript{r}, but also have given an algorithm to find the exact positions of these combinations in the matrix of all possible combinations(\textit{prast\=ara}) of n G-Ls. 

\section{Conclusion}
	The use of mathematical algorithms and of recursion dates back to around 200 B.C..  \textit{Pi\.nga\d{l}a}  used recursion extensively to describe the algorithms. Further, the use of stack to store the information of intermediate operations, in \textit{Pi\.nga\d{l}a}'s algorithms is also worth mentioning. All these algorithms use a terminating condition also, ensuring that the recursion terminates. Recursive algorithms are easy to conceptualise, and implement mechanically. 
We notice the use of method of recursion and also the binomial expansion in the later works on mathematics such as \textit{brahmasphotasiddhanta}\superscript{7} with commentary by \textit{pruthudaka} on summing a geometric series, or \textit{Bhattotpala}'s commentary on \textit{bruhatsamhit\=a} etc. 
They present a mathematical model corresponding to the algorithm. However,  the iterative algorithms are easy from user's point of view. They are directly executable for a given value of inputs, without requirement of any stacking of variables. Hence the later commentators such as \textit{Ked\=ar} \textit{Bhatt} might have used only iterative algorithms.  The \textit{s\=utra} style was prevalent in India, and unlike modern mathematics, the Indian mathematics was passed from generations to generations verbally, through \textit{s\=utra}s. Sutras being very brief, and compact, were easy to memorise and also communicate orally. However, it is still unexplored what features of Natural language like Sanskrit have Indian mathematicians used for mathematics as opposed to a specially designed language of modern mathematics that made the Indian mathematicians communicate mathematics orally effortlessly.

\section{Acknowledgment}
	The material in this paper was evolved while teaching a course on “Glimpses of Indian Mathematics” to the first year students of the Integrated Masters course at the University of Hyderabad. The author acknowledges the counsel of her father A. B. Padmanabha Rao, Subhash Kak, and G\`erard Huet who gave useful pointers and suggestions.

{}
\end{document}